\documentclass[a4paper,12pt]{article}

\usepackage{amsmath,amssymb}
\usepackage{graphicx} 
\usepackage{epsfig}
\usepackage{amsfonts}
\usepackage{amstext}

\newcommand{\beq}{\begin{equation}}
\newcommand{\eeq}{\end{equation}}
\newcommand{\done}{\hspace*{\fill} $\Box$}
\newcommand{\RE}{I\!\!R}

\newtheorem{theorem}{Theorem}
\newtheorem{remark}{Remark}
\newtheorem{corollary}{Corollary}

\newtheorem{lemma}{Lemma}
\newtheorem{example}{Example}

\title{\bf An Inverse Optimality Method to Solve a Class of Optimal Control Problems}

\begin{document}

\author{Luis Rodrigues$^1$, Didier Henrion$^{2,3,4}$, Mehdi Abedinpour Fallah$^1$}

\footnotetext[1]{Department of Electrical and Computer Engineering, Concordia University,
1515 St. Catherine Street, EV12.111, Montr\'eal, QC H3G 2W1, Canada.}
\footnotetext[2]{CNRS; LAAS; 7 avenue du colonel Roche, F-31077 Toulouse; France.}
\footnotetext[3]{Universit\'e de Toulouse; UPS, INSA, INP, ISAE; UT1, UTM, LAAS; F-31077 Toulouse; France}
\footnotetext[4]{Faculty of Electrical Engineering, Czech Technical University in Prague,
Technick\'a 2, CZ-16626 Prague, Czech Republic}
\maketitle

\begin{abstract}
This paper presents an inverse optimality method to solve the Hamilton-Jacobi-Bellman equation for a class of nonlinear problems for which the cost is quadratic and the dynamics are affine in the input.
The method is inverse optimal because the running cost that renders the control input optimal is also explicitly determined.
One special feature of this work, as compared to other methods in the literature, is the fact that the solution is obtained directly for the control input. 
The value function can also be obtained after one solves for the control input.
Furthermore, a Lyapunov function that proves at least local stability of the controller is also obtained.
In this regard the main contribution of this paper can be interpreted in two different ways: offering an analytical expression for Lyapunov functions for a class of nonlinear systems and obtaining an optimal controller for the same class of systems using a specific optimization functional. 
We also believe that an additional contribution of this paper is to identify explicit classes of systems and optimization functionals for which optimal control problems can be solved analytically.
In particular, for second order systems three cases are identified: i) control input only as a function of the second state variable, ii) control input affine in the second state variable when the dynamics are affine in that variable and iii) control input affine in the first state variable when the dyamics are affine in that variable.  
The relevance of the proposed methodology is illustrated in several examples, including the Van der Pol oscillator, mass-spring systems and vehicle path following.
\end{abstract}

\begin{center}
Hamilton-Jacobi-Bellman equation.
\end{center}

\section{Introduction}
\noindent Optimal control problems are hard to solve because the optimal controller is the solution of a partial differential equation called the Hamilton-Jacobi-Bellman (HJB) equation \cite{lqr}.
However, when the cost is quadratic and the dynamics are affine in the input there is an explicit solution for the input as a function of the derivatives of the value function.
This fact will be used to develop a method to solve the HJB equation for a class of nonlinear systems.
The main motivation for this work comes from the controller designer's perspective.
When designers are faced with a control engineering problem and want to formulate it in the optimal control framework, the first challenge is to choose the most appropriate cost that will yield a control solution with physical significance. Although this is a difficult choice, quite often the following three properties are required for the design:
\begin{enumerate}
\item The closed loop system should be asymptotically stable to a desired equilibrium point
\item The system should have enough damping so that the trajectories do not take too long to settle around the desired equilibrium point 
\item The control energy should be penalized in the cost to avoid high control inputs that can saturate actuators
\end{enumerate} 
The particular functions involved in the cost are not usually pre-defined, except possibly the requirement on the control energy that is usually represented by a quadratic cost on the input.
The work on this paper attempts to find a controller and a cost that together meet the requirements 1--3 and render the controller optimal relative to that cost.
To that aim, the cost will be fixed to be quadratic in the input and have an unknown term in the state that shall be determined.
The solution is therefore based on the concept of inverse optimality.
One special feature of this method, as compared to other methods in the literature, is the fact that the solution is obtained directly for the control input without needing to assume or compute a value function first. Rather, the value function is obtained after one has solved for the control input. A Lyapunov function will also be constructed, at least locally.
Work on optimal control and approximate solutions, such as inverse optimality, has started in the sixties (see for example \cite{kalman}, \cite{lukes} and references therein), concentrating mostly on linear quadratic problems driven by aerospace applications. Thirty years later, the concept of inverse optimality has been revisited by many authors to address nonlinear optimal control problems.
In a pioneering paper, Lukes \cite{lukes} approximates the solution to an optimal control problem with analytic functions by a a Taylor series, starting with first order terms in the dynamics and second order terms in the cost. The resulting controller is therefore the sum of a Linear Quadratic Regulator (LQR) with higher order terms.
Reference \cite{doyle} finds the dynamics that verify the HJB equation given the running cost and a value function. In \cite{kokotovic} an analytical expression for a stabilizing controller is obtained for feedback linearizable dynamics given the coordinate transformation that feedback linearizes the system, a control Lyapunov function obtained as the solution of the Riccatti equation for the linearized dynamics and a bound on the decay rate of the Lyapunov function. It is shown that the controller is optimal relative to a cost involving a control penalty bias. Reference \cite{margaliot} uses Young's inequality, which was used before in \cite{krstic} for the design of input-to-state stabilizing controllers, to find an analytical expression for the solution to a class of nonlinear optimal control problems.
An expression for the cost that makes the controller optimal was also found. However, there is no indication as to what conditions must be satisfied such that the obtained cost is a sensible cost, namely, such that it is non-negative. This is shown on a case-by-case basis in the examples.
Reference \cite{krstic} showed that both the inverse optimal gain assignment and ${\mathcal H}_{\infty}$ problems are solvable for the case where the system is in strict-feedback nonlinear form. For a similar strict-feedback nonlinear form, the work presented in \cite{kokotovic2} develops a recursive backstepping controller design procedure and the corresponding construction of the cost functional using nonlinear Cholesky factorization.
It is shown that under the assumptions that the value function for the system has a Cholesky factorization and the running cost is convex, it is possible to construct globally stabilizing control laws to match the optimal ${\mathcal H}_{\infty}$ control law up to any desired order, and to be inverse optimal with respect to some computable cost function.
In terms of applications, reference \cite{rigidspacecraft} presents an inverse optimal control approach for regulation of a rotating rigid spacecraft by solving an HJB equation. The resulting design includes a penalty on the angular velocity, angular position, and the control torque. The weight in the penalty on the control depends on the current state and decreases for states away from the origin. Inverse optimal stabilization of a class of nonlinear systems is also investigated in \cite{26thChinese} resulting in a controller that is optimal with respect to a meaningful cost function. The inverse optimality approach used in \cite{rigidspacecraft} and \cite{26thChinese} requires the knowledge of a control Lyapunov function and a stabilizing control law of a particular form. In \cite{pendulum} an optimal feedback controller for bilinear systems is designed to minimize a quadratic cost function. This inverse optimal control design is also applied to the problem of the stabiliz!
 ation of an inverted pendulum on a cart with horizontal and vertical movement.

Building on the concept of inverse optimality, but in contrast with previous approaches, the objective of this paper is to offer a solution method for a class of nonlinear systems that can determine at the same time a controller and a sensible non-negative cost that renders the controller optimal. Although limited to models up to third order, an important contribution of the work presented in this paper is the fact that the models considered here do not have to be in strict nonlinear feedback form. In fact, the derivative of state variable $i$ does not necessarily have to be an affine function of state variable $i+1$ for the models considered in this paper. Furthermore, the running cost is not assumed to be convex. In addition, the analytical solution for the control input is obtained directly, without needing to first assume or compute any coordinate transformation, value function, or Lyapunov function. The value function and a Lyapunov function can however be computed once!
  the optimal control input has been found.
Finally, conditions are given such that the cost that makes the controller optimal is a sensible non-negative cost.
The paper is organized as follows. First the optimal control problem is defined and solved for a class of second order systems, followed by its extension to a class of third order systems and conclusions. Several examples are presented throughout the paper.
In the notation used in the paper $V_{x_i}$ denotes the partial derivative of $V$ with respect to $x_i$ and $f'(x_i)$ denotes the derivative of function $f$ with respect to its only argument $x_i$.

\section{Optimal Control Problem Definition and Solution: Second Order Systems}\label{twostates}
\noindent Consider the following optimal control problem
\beq
\begin{array}{rclr}
V(x_0)& = & \inf \int_0^{\infty} \left\{Q(x)+ru^2\right\}dt&\\
& s.t.& \dot x_1(t) = f_1(x_1,x_2)&\label{nlsystemeq4}\\
			&& \dot x_2(t) = f_2(x_1,x_2)+bu&\\
			&& x(0)=x_0,~u\in{\mathcal U}&
\end{array}
\eeq
where it is assumed that $V$ is of class ${\mathcal C}^1,~Q(x)\ge 0,~Q(0)=0,~r>0,~b\neq 0,~x(t)=[x_1(t)~x_2(t)]^T\in\RE^2$.
The set $\mathcal U$ represents the allowable inputs, which are considered to be Lebesgue integrable functions. 
The functions $f_1,~f_2$ are not identically zero and are assumed to be continuous with $f_1(0)=f_2(0)=0$.
These functions will be further constrained in the theorems presented in the paper.
The term
\beq\label{running_cost}
L(x_1,x_2,u)=Q(x)+ru^2
\eeq
is called the running cost.
When $f_1,~f_2$ are linear, from the LQR theory \cite{lqr}, one knows that the optimal solution is a linear state feedback law $u=-k_1x_1-k_2x_2$.
Inspired by this fact, for nonlinear $f_1,~f_2$ we will search for nonlinear additive state feedback laws of the form $u=u_1(x_1)+u_2(x_2)$ with $u_1(0)=u_2(0)=0$. The first problem to be solved is to find out for what forms of $Q(x)$ such a control input exists. The second problem is to find a solution $u=u_1(x)+u_2(x)$ given a $Q(x)$ in the allowed form.
We start by presenting necessary conditions that the value function $V$ must verify for additive control solutions to exist.
\newline
\begin{lemma}\label{necessary}
Assume that a control solution of the form 
\begin{equation}\label{control}
u(x)=u_1(x_1)+u_2(x_2)
\end{equation}
with $u_1(x_1)$ of class ${\mathcal C}^1$ and $u_2(x_2)$ continuous, exists for problem (\ref{nlsystemeq4}) with $u_1(0)=u_2(0)=0$. Furthermore, assume that a class ${\mathcal C}^1$ function $V$
exists that verifies the corresponding HJB equation
\beq\label{hjbvector}
\inf_u H\left(x_1,x_2,u,V_{x_1},V_{x_2}\right)=0
\eeq
where
\beq\label{Hamiltonian}
\begin{array}{rcl}
H&=&Q(x)+ru^2+V_{x_1}f_1(x_1,x_2)+V_{x_2}f_2(x_1,x_2)+V_{x_2}bu
\end{array}
\eeq
with boundary condition $V(0)=0$. 
Then $V$ must be of the form
\beq\label{valuefunction}
V(x)=-2b^{-1}r\left(x_2u_1(x_1)+U_2(x_2)\right)+h(x_1)
\eeq
where $u_1(x_1),~h(x_1)$ and $U(x_2)$ are functions of class ${\mathcal C}^1$ with
\begin{equation}
u_2(x_2) = U_2'(x_2),\quad h(0)= 2b^{-1}rU_2(0)\label{necessarycond1b}
\end{equation} 
Furthermore, $u_1$ and $u_2$ are solutions of the equation
\beq\label{hjbnecessary}
\begin{array}{rcl}
Q-ru_1^2-2ru_1u_2-ru_2^2-2b^{-1}rx_2u_1'f_1+h'f_1-2b^{-1}ru_1f_2-2b^{-1}ru_2f_2 &=& 0
\end{array}
\eeq
where the arguments of the functions were omitted for simplicity.
\end{lemma}
{\bf Proof.}
Consider the HJB equation (\ref{hjbvector}) associated with (\ref{nlsystemeq4}).
The necessary condition on $u$ to be a minimizer is
\beq\label{dvdx2}
V_{x_2}=-2b^{-1}ru(x)
\eeq
and therefore
\beq\label{valuetwostates}
V(x)=-2b^{-1}r\int u(x)dx_2 + h(x_1)
\eeq
where $h(x_1)$ is an arbitrary integration function of $x_1$.
Replacing (\ref{control}) into (\ref{valuetwostates}) yields (\ref{valuefunction}).
From the boundary condition $V(0)=0$ one obtains the constraint (\ref{necessarycond1b}) taking into account that $u_1(0)=0$.
Differentiating (\ref{valuetwostates}) with respect to $x_1$ and using (\ref{control}) yields
\beq\label{dvdx1}
V_{x_1}=-2b^{-1}rx_2u_1'(x_1)+h'(x_1)
\eeq
Finally, replacing (\ref{control}), (\ref{dvdx2}), and (\ref{dvdx1}) in (\ref{hjbvector}) yields (\ref{hjbnecessary}) after rearranging.
This finishes the proof. \done

\begin{remark}
It is important to note that assuming a control input of the form (\ref{control}) allows one to transform the HJB equation into an ordinary differential equation instead of a partial differential equation. Furthermore, it is interesting to note that if the value function (\ref{valuefunction}) does not have cross terms in $x_1$ and $x_2$, from (\ref{dvdx2}), the controller will only depend on $x_2$.
\end{remark}
\vspace{10pt}
Based on the form of (\ref{hjbnecessary}), this equation will now be solved for three different cases: 
i) control input only as a function of $x_2$, ii) control input affine in $x_2$ when the dynamics are affine in that variable and iii) control input affine in $x_1$ when the dyamics are affine in that variable.  

\subsection{Case I: Solutions depending only on $x_2$} 
For this case we first assume that $f_2$ is only a function of $x_2$.
The result is stated in Theorem \ref{case1}.
\newline
\begin{theorem}\label{case1}
Assume that $f_1(x_1,x_2)$ and $f_2(x_2)$ are continuous and such that
\beq\label{plantconstraint}
f_1(0,0)=0,~f_2(0)=0
\eeq
and $f_1$ is not identically zero. If $Q(0,0)=0$ and $Q$ is of the form
\begin{equation}\label{Qcase1}
Q(x_1,x_2)=-g(x_1)f_1(x_1,x_2)+Q_2(x_2)
\end{equation}
where $g$ is a function of class ${\mathcal C}^1$ not identically zero, $Q_2\neq 0$ and
\begin{eqnarray}\label{Qconditions}
Q_2(x_2)&\ge& 0\nonumber\\
-g(x_1)f_1(x_1,x_2)+Q_2(x_2)&\ge& 0
\end{eqnarray}
then the stabilizing control input $u=u_2(x_2)$ that is a solution of the quadratic equation
\begin{equation}\label{Q2case1}
Q_2(x_2)-ru_2^2-2b^{-1}ru_2f_2(x_2)=0
\end{equation}
is an optimal solution of problem (\ref{nlsystemeq4}) if it is continuous and the corresponding value function is given by
\beq\label{valuecase1}
V(x_1,x_2)=-2b^{-1}r\int u_2(x_2)dx_2 + \int g(x_1)dx_1 - 2rb^{-1}U_2(0)
\eeq
Furthermore, if $u_2$ is of class ${\mathcal C}^1$ and
\begin{equation}\label{posdefcase1}
g'(x_1)> 0,~x_1\neq 0\quad u_2'(x_2) < 0,~x_2\neq 0
\end{equation}
then $V$ is positive definite and it is a local Lyapunov function.
The function $V$ is a global Lyapunov function if it is radially unbounded.
Finally, the trajectories converge to one of the minimizers of $L(x_1,x_2,u(x_1,x_2))$, i.e, to a point $(x_1,x_2)$ such that $L=0$. If $L$ is convex, then the trajectories will converge to the origin for all initial conditions.
\end{theorem}
{\bf Proof.}
From the proof of Lemma \ref{necessary} the HJB equation can be written as (\ref{hjbnecessary}).
With $u_1=0$ this equation becomes
\beq\label{hjbvector3}
\begin{array}{rcl}
Q-ru_2^2+h'f_1-2b^{-1}ru_2f_2 &=& 0
\end{array}
\eeq
where $Q\ge 0$ under conditions (\ref{Qconditions}).
Making
\beq\label{hprimecase1}
h'(x_1)=g(x_1)
\eeq
using (\ref{Qcase1}) and (\ref{Q2case1}) yields $0=0$, and therefore the HJB equation is satisfied.
The HJB equation is a sufficient condition for the control input (\ref{control}) with $u_1(x_1)=0$ to be a solution that minimizes the cost of problem (\ref{nlsystemeq4}) because the second derivative of the Hamiltonian (\ref{Hamiltonian}) with respect to $u$ is equal to $2r>0$.
Using $u_1=0$ and replacing the integral of (\ref{hprimecase1}) in (\ref{valuefunction}) yields the value function (\ref{valuecase1}) taking into account (\ref{necessarycond1b}). 
Observe that from the HJB equation (\ref{hjbvector}) and from $Q(x)\ge 0$, if $u_2$ is continuous we have
\beq\label{dotlyap}
\dot V = -L(x_1,x_2,u)\le 0
\eeq
which makes $V$ a local Lyapunov function for the system if $u_2$ is also of class ${\mathcal C}^1$ because of the conditions (\ref{posdefcase1}) on the Hessian of $V$.
If $V$ is also radially unbounded it is a global Lyapunov function.
Finally, since the optimal cost (\ref{valuecase1}) is finite for all initial conditions, then the trajectories will converge to one of the minimizers of $L(x_1,x_2,u(x_1,x_2))$ because $L\ge 0$ and $\lim_{t\to\infty}L=0$ for integrability. If $L$ is convex, then the trajectories must converge to the origin because the origin is the only minimizer of $L$.
This finishes the proof. \done
\newline

\begin{remark}
Note that equation (\ref{Q2case1}) with $Q_2(x_2)\ge 0$ corresponds to the solution of an optimal control problem with running cost $L=Q_2(x_2)+ru^2$ and first order dynamics $\dot x_2=f_2(x_2)+bu_2$. Therefore, the result of Theorem \ref{case1} reduces the solution of an optimal control problem for a second order system to the solution of an optimal control problem for a first order system.
\end{remark}

\begin{example}
If $f_1(x_1,x_2)=-x_1^3-2x_1x_2,~f_2(x_2)=x_2\sqrt{3\left(1+x_2^2\right)}$ and $Q(x_1,x_2)=\left(x_1^2+x_2\right)^2+x_2^4,~b=r=1$ then using the result of Theorem \ref{case1} we get $g(x_1)=x_1,~Q_2(x_2)=x_2^2+x_2^4$ and $u=u_2(x_2)=-\left(2+\sqrt{3}\right)x_2\sqrt{1+x_2^2}$.
\end{example}

We now assume that $f_1$ is only a function of $x_2$.

\begin{theorem}\label{case1b}
Assume that $f_1(x_2)$ and $f_2(x_1,x_2)=f_{21}(x_1)+f_{22}(x_1,x_2)$ are continuous and such that
\beq\label{plantconstraint1b}
f_1(0)=0,~f_{22}(0,0)=0
\eeq
and $f_1,~f_{21}$ are not identically zero.
If $Q$ is of the form
\begin{equation}\label{Qcase1b}
Q(x_1,x_2)=rk^2f_1^2+2b^{-1}rkf_1f_{22}
\end{equation}
and
\begin{eqnarray}\label{fconditions1b}
b^{-1}kf_1(x_2)f_{22}(x_1,x_2)&\ge& 0
\end{eqnarray}
then the control input $u=u_2(x_2)=kf_1(x_2)$  is an optimal solution of problem (\ref{nlsystemeq4}) and the corresponding value function is given by
\beq\label{valuecase1b}
V(x_1,x_2)=-2b^{-1}rk\left(\int f_1(x_2)dx_2 - \int f_{21}(x_1)dx_1\right)
\eeq
Furthermore, if $f_1,~f_{21}$ are of class ${\mathcal C}^1$ and
\begin{eqnarray}\label{posdefcase1b}
b^{-1}kf_1'(x_2) &<& 0,~x_2\neq 0\nonumber\\
b^{-1}kf_{21}'(x_1) &>& 0,~x_1\neq 0
\end{eqnarray}
then $V$ is positive definite and it is a local Lyapunov function.
The function $V$ is a global Lyapunov function if it is radially unbounded.
Finally, the trajectories converge to one of the minimizers of $L$. If $L$ is convex, then the trajectories will converge to the origin for all initial conditions.
\end{theorem}
{\bf Proof.}
From the proof of Theorem \ref{case1} with $u_1=0$ the HJB equation can be written as
\beq\label{hjbvector31b}
\begin{array}{rcl}
Q-ru_2^2(x_2)+h'(x_1)f_1(x_2)-2b^{-1}ru_2(x_2)\left[f_{21}(x_1)+f_{22}(x_1,x_2)\right] &=& 0
\end{array}
\eeq
Making $u_2(x_2)=kf_1(x_2)$,
\beq\label{hprimecase1b}
h'(x_1)=2b^{-1}rkf_{21}
\eeq
and using (\ref{Qcase1b}) yields $0=0$, and therefore the HJB equation is satisfied.
Note that under assumption (\ref{fconditions1b}), the running cost $L$ is non-negative.
The rest of the proof follows the same reasoning of the proof of Theorem \ref{case1}.\done

\begin{example}
If $f_1(x_2)=x_2^3,~f_2(x_1,x_2)=-x_1^3-x_1^2x_2,~Q(x_1,x_2)=rk^2x_2^6-2b^{-1}rkx_1^2x_2^4$ then using the result of Theorem \ref{case1b} we get that $u=kx_2^3$ is the optimal control with $V(x_1,x_2)=-2b^{-1}rk\left(x_2^4/4 + x_1^4/4\right)$ where $k$ is chosen such that $b^{-1}k<0$.
\end{example}

\subsection{Case II: Solutions that are affine in $x_2$ and depend on both $x_1,x_2$}
For this case we assume that both $f_1$ and $f_2$ are affine functions of $x_2$.
The main result is stated in Theorem \ref{case2}.
\begin{theorem}\label{case2}
Assume that
\begin{eqnarray}
f_1(x_1,x_2)&=&g_1(x_1)+g_2(x_1)x_2\nonumber\\
f_2(x_1,x_2)&=&g_3(x_1)+g_4(x_1)x_2\nonumber\\
\end{eqnarray}
where $g_2(x_1)\neq 0$, $g_3(x_1)$, $g_4(x_1)$ are continuous functions, $g_1(x_1)$ is of class ${\mathcal C}^1,~g_1(0)=g_2(0)=g_3(0)=g_4(0)=0$.
If given $Q_1(x_1)\ge0,q_2>0$, the stabilizing solution $u_1$ of
\begin{equation}\label{usolution1}
Q_1(x_1)-ru_1^2-2b^{-1}ru_1g_3(x_1)=0
\end{equation}
is of class ${\mathcal C}^1$ then the control input
\beq\label{equationu1}
u=u_1(x_1)-k_2x_2
\eeq
with $k_2=\pm\sqrt{q_2r^{-1}},~b^{-1}k_2>0$ is a solution of the optimal control problem (\ref{nlsystemeq4}) when $Q$ is of the form 
\beq\label{Qcase2}
Q(x)=Q_1(x_1)+q_2x_2^2+2rb^{-1}\left(u'_1g_2-k_2g_4\right)x_2^2-h'g_1
\eeq
and
\beq\label{k22}
k_2^2+2b^{-1}\left(u'_1g_2-k_2g_4\right)\ge 0,~h'g_1\le 0.
\eeq
where $h(x_1)$ is a function of class ${\mathcal C}^1$ satisfying
\beq\label{h}
h'g_2=-2rk_2\left(u_1+b^{-1}g_3\right)+2rb^{-1}\left(u_1g_4+u'_1g_1\right)
\eeq
The resulting value function is
\beq\label{finalvaluecase2}
V(x)=rb^{-1}\left(-2u_1x_2+k_2x_2^2\right) - 2rk_2 \int g^{-1}_2\left(u_1+b^{-1}g_3\right)dx_1 + 2rb^{-1}\int g^{-1}_2\left(u_1g_4+u'_1g_1\right)dx_1+c
\eeq
where $c$ is chosen such that the boundary condition
$V(0)=0$ is satisfied.
The function $V$ is also a local Lyapunov function provided it is positive definite in a region
around the origin. If $V$ is globally positive definite and radially unbounded then it is a Lyapunov
function.
If $L$ is convex, then the trajectories will converge to the origin for all initial conditions.
If $L$ is not convex then the trajectories will converge to one of the minimizers of $L$.
\end{theorem}
\vspace{10pt}

\noindent{\bf Proof.}
Taking into account (\ref{usolution1}), (\ref{equationu1}) and (\ref{Qcase2}), equation (\ref{hjbnecessary}) becomes after rearranging
\beq\label{hjbvector4}
q-r\left(k_2^2+2b^{-1}\left(u'_1g_2-k_2g_4\right)-\frac{q_2}{r}\right)x_2^2+2rx_2\left(k_2\left(u_1+b^{-1}g_3\right)-b^{-1}\left(u_1g_4+u'_1g_1\right)+\frac{h'g_2}{2r}\right)+h'g_1=0.
\eeq
where $q=2rb^{-1}\left(u'_1g_2-k_2g_4\right)x_2^2-h'g_1$. Replacing $k_2=\pm\sqrt{q_2r^{-1}}$ and $h'$ given by (\ref{h}) into (\ref{hjbvector4}), the HJB equation is satisfied.
Note that this is a sufficient condition for optimality because the second derivative of the Hamiltonian with respect to $u$ is equal to $2r>0$. 
For positivity of the running cost $L$ one must have
\beq\label{k2}
Q_1(x_1)+q_2x_2^2+2rb^{-1}\left(u'_1g_2-k_2g_4\right)x_2^2-h'g_1\ge 0.
\eeq
Note that this constraint is always satisfied if (\ref{k22}) holds.
The value function $V$ is obtained replacing the control inputs and $h$ in (\ref{valuefunction}).
Note that from the Hessian of $V$, $b^{-1}k_2$ is one of the sufficient conditions for $V$ to be positive definite.
The rest of the proof follows the same argument as in the proof of Theorem \ref{case1}.
\done

\begin{remark}
It is important to note that the equation (\ref{usolution1}) corresponds to the solution of an optimal control problem with running cost $Q_1(x_1)$ and first order dynamics $\dot x_1=g_3(x_1)+bu_1$. Therefore, the result of Theorem \ref{case2} reduces the solution of an optimal control problem for a second order system to the solution of an optimal control problem for a first order system plus the addition of a viscous damping term $u_2(x_2)=-k_2x_2$.
\end{remark}

\begin{remark}
It is interesting to note that when $g_1(x_1)=0,~g_2(x_1)=1,~Q_1(x_1)=0$ and $x_1g_3(x_1)< 0,~x_1\neq 0$, meaning that $\dot x_1 = g_3(x_1)$ is asymptotically stable, then the result of Theorem \ref{case2} coincides with the result of Theorem \ref{case1b}.
\end{remark}

\begin{example}
Consider the mass-spring system with dynamics $g_1(x_1)=0$, $g_2(x_1)=1$, $g_3(x_1)=-x_1^3$, $g_4(x_1)=0$, $b=1$, and assume $Q_1(x_1)=0$. 
Then, using the results of Theorem \ref{case2}, from (\ref{usolution1}) we get $u_1=0$.
Therefore, with the running cost $L(x,u)=q_2x_2^2+ru^2$ the solution is $u=-\sqrt{q_2r^{-1}}x_2$ with value function $V(x)=\sqrt{q_2r}\left(x_2^2+0.5x_1^4\right)$, which is also a Lyapunov function for the closed loop system. Note that the control input is adding viscous damping to the mass-spring system to stabilize it to the origin, which makes perfect sense from a physical point of view.
\end{example}

\begin{example}
Consider the Van der Pol oscillator with dynamics given by $b=1,~g_1(x_1)=0,~g_2(x_1)=1,~g_3(x_1)=-x_1,~g_4(x_1)=0.5(1-x_1^2)$, and assume $Q_1(x_1)=0,~q_2=1,~r=1$. Then, using the results of Theorem \ref{case2}, from (\ref{usolution1}) we get $u_1=0$ and the optimal controller $u=-x_2$ with associated value function $V(x)=x_1^2 + x_2^2$, which is also a Lyapunov function for the closed loop system. The running cost is $L(x,u)=x_1^2x_2^2+u^2$. This controller makes perfect sense from a physical point of view because to damp out the oscillations and make the trajectories converge to the origin the input simply adds viscous damping.
\end{example}

\begin{example}
Let $b=1,~g_1(x_1)=-x_1^3,~g_2(x_1)=1,~g_3(x_1)=g_4(x_1)=0,~Q_1(x_1)=q_1x_1^2,~r=1,~q_1,~q_2>0$.
This is a system in strict feedback form to which backstepping techniques can be applied. 
From the results of Theorem \ref{case2}, solving (\ref{usolution1}) the resulting controller is $u=-\sqrt{q_1}x_1-\sqrt{q_2}x_2$. From (\ref{h}) one gets 
$$
h'g_1=-2\left(\sqrt{q_1q_2}x_1^4+\sqrt{q_1}x_1^6\right)\le 0
$$ 
and from (\ref{Qcase2})
$$
Q(x_1,x_2)=q_1x_1^2+2\sqrt{q_1q_2}x_1^4+2\sqrt{q_1}x_1^6+\left(q_2-2\sqrt{q_1}\right)x_2^2
$$
The constraint (\ref{k22}) is $q_2\ge 2\sqrt{q_1}$.
Finally, from (\ref{finalvaluecase2}), the value function is
$$
V(x)=2\sqrt{q_1}x_1x_2+\sqrt{q_2}\left(x_2^2+\sqrt{q_1}x_1^2\right)+2\sqrt{q_1}\frac{x_1^4}{4}
$$
which is a Lyapunov function for the closed loop system.
\end{example}

\subsection{Case III: Solutions that are affine in $x_1$ and depend on both $x_1,~x_2$}
For this case we assume that both $f_1$ and $f_2$ are affine in $x_1$.
The main result of this section is stated in the next theorem.
\newline
\begin{theorem}\label{mainresult}
Assume that there exist real scalars $a,~b,~c,~d$ such that
\begin{eqnarray}
f_1(x_1,x_2)&=&ax_1+f(x_2)\nonumber\\
f_2(x_1,x_2)&=&cx_1+df(x_2)
\end{eqnarray}
where $f$ is not identically zero and is assumed to be continuous with $f(0)=0$ and with a locally positive definite anti-derivative $F(x_2)$ such that $F'(x_2)=f(x_2)$. Assume further that $c(ad-c)\ge 0$ and that for some $\beta>0$
\beq\label{plantconstraint3}
\beta a^2\ge c^2 \ge a^2d^2
\eeq
This implies that either $a=c=0$ or $a\neq 0, c\neq 0$ or $a\neq 0, c=0, d=0$.
Furthermore, assume that
\beq
Q(x)=q_1x_1^2+q_2x_2^2+q(x)
\eeq
where $q_1\ge 0,~q_2>0$ and
\beq\label{q1}
\left\{
\begin{array}{lr}
q_1=q_2c^2a^{-2}+2rb^{-2}c\left(ad-c\right), & a\neq 0\\
q_1\ge q_2d^2, & a=0 
\end{array}
\right.
\eeq
Finally, let $q(x)$ be chosen as
\beq\label{Q}
q=2rk_1k_2x_1x_2+rb^{-2}\left(k_1^2k_2^{-2}-d^2\right)f^2
\eeq
Then, there exist gains $k_1,~k_2,~k$ verifying
\beq\label{k2constraints}
k_2=\sqrt{\frac{q_2}{r}}
\eeq
\beq\label{k1constraints}
k_1=\left\{
\begin{array}{lr}
-ca^{-1}k_2, & a\neq 0\\ 
\sqrt{\frac{q_1}{r}}, & a=0 
\end{array}
\right.
\eeq
\beq\label{kconstraint}
k=b^{-1}\left(d+\frac{k_1}{k_2}\right)
\eeq
such that the control input (\ref{control}) is a solution of the HJB equation (\ref{hjbvector}) associated with (\ref{nlsystemeq4}) with value function
\beq\label{finalvalue}
\begin{array}{rcl}
V(x)&=&-rb^{-1}kcx_1^2+rb^{-1}\left(\frac{k_1}{\sqrt{k_2}}x_1+\sqrt{k_2}x_2\right)^2\\
&&+2rb^{-1}k\left(F(x_2)-F(0)\right)
\end{array}
\eeq
The function $V(x)$ is also a local Lyapunov function provided $b>0,~kc\le 0$ and the term $kF(x_2)$ is locally positive definite, i.e, if for some class ${\mathcal K}$ function $\delta$ and positive $\gamma$ one has
\begin{equation}\label{convexity}
kF(x_2)\ge\delta\left(\|x_2\|\right),~\forall x_2\in\Omega\triangleq\{x_2~:~\|x_2\|\le\gamma\}
\eeq
If $V$ is globally positive definite and radially unbounded then it is a Lyapunov function.
Finally, the trajectories converge to one of the minimizers of $L$. If $L$ is convex, then the trajectories will converge to the origin for all initial conditions.
\end{theorem}
{\bf Proof.}
The HJB equation (\ref{hjbnecessary}) can be written as
\beq\label{hjbvector2}
\begin{array}{rcl}
0&=&\left(q_1-rp_1\right)x_1^2+\left(q_2-rk_2^2\right)x_2^2+q\\
&&+h'f+p_2x_1f-rk\left(k-2b^{-1}d\right)f^2\\
&&-2r\left[k_1k_2-b^{-1}\left(k_1a+k_2c\right)\right]x_1x_2+h'ax_1\\
&&+2r\left[b^{-1}\left(k_1+k_2d\right)-k_2k\right]x_2f
\end{array}
\eeq
where the arguments of the functions were omitted for simplicity and
\beq\label{ps}
\begin{array}{rcl}
p_1(k_1,k)&=&k_1^2-2b^{-1}ck_1\\
p_2(k_1,k)&=&2r\left[b^{-1}\left(k_1d+kc\right)-k_1k\right]
\end{array}
\eeq
Since by assumption $q_2>0$, then (\ref{k2constraints}) implies $k_2\neq 0$ and (\ref{kconstraint}) is well defined.
Note that the term $\left[b^{-1}\left(k_1+k_2d\right)-k_2k\right]x_2f$ in (\ref{hjbvector2}) vanishes because of (\ref{kconstraint}).
Note that if $a=0$ then $c=0$ because of inequalities (\ref{plantconstraint3}).
This observation together with (\ref{k1constraints}) yields
$$
k_1a+k_2c=0
$$
and therefore the term $\left(k_1a+k_2c\right)x_1x_2$ in (\ref{hjbvector2}) vanishes.
Making
\beq\label{hprime}
h'(x_1)=-p_2(k_1,k)x_1
\eeq
the term $h'f+p_2x_1f$ in (\ref{hjbvector2}) vanishes.
Using (\ref{ps}), (\ref{q1}), (\ref{k2constraints}), (\ref{k1constraints}), and (\ref{kconstraint}) for the case $a\neq 0$, and using (\ref{ps}) and (\ref{k1constraints}) for the case $a=0$ (which implies also $c=0$), one finds that $q_1-rp_1=ap_2$.
We also see that the term $\left(q_1-rp_1\right)x_1^2+h'ax_1$ in (\ref{hjbvector2}) vanishes.
The term $\left(q_2-rk_2^2\right)x_2^2$ vanishes because of constraint (\ref{k2constraints}).
Using (\ref{Q}) and (\ref{kconstraint}) the term $q-2rk_1k_2x_1x_2-rk\left(k-2b^{-1}d\right)f^2$ in (\ref{hjbvector2}) also vanishes.
Since all terms in (\ref{hjbvector2}) vanish, the HJB equation is satisfied.
This is a sufficient condition for the control input (\ref{control}) to be a solution that minimizes the cost of problem (\ref{nlsystemeq4}) because the second derivative of the Hamiltonian (\ref{Hamiltonian}) with respect to $u$ is equal to $2r>0$.
The running cost is a sensible cost because from (\ref{running_cost}) and (\ref{plantconstraint3})--(\ref{k1constraints}) it is given by
\begin{equation}\label{running}
L=r\left(k_1x_1+k_2x_2\right)^2+2rc(ad-c)b^{-2}x_1^2+rb^{-2}\left(k_1^2k_2^{-2}-d^2\right)f^2+ru^2\nonumber
\end{equation}
and it is non-negative with a minimun at $x_1=x_2=u=0$ under the assumptions $~c(ad-c)\ge 0$, (\ref{plantconstraint3}), (\ref{q1}), (\ref{k2constraints}), (\ref{k1constraints}).
Replacing the integral of (\ref{hprime}) in (\ref{valuefunction}), using (\ref{ps}) and (\ref{kconstraint}) yields the value function (\ref{finalvalue}). 
The boundary condition $V(0)=0$ yields the term $-2rb^{-1}kF(0)$, which is a constant of integration. 
The rest of the proof follows the same argument as in the proof of Theorem \ref{case1}. \done
\newline

\begin{remark}
It is interesting that the square of the nonlinearity comes naturally as a term in the cost, although this would be difficult to predict based on a general tendency to always construct costs that have only quadratic terms on the state.
\end{remark}
\vspace{10pt}

\begin{example}
For system (\ref{nlsystemeq4}) with $f(x_2)=x_2,~a=c=d=0$ and $b=1$ one obtains a double integrator.
According to Theorem \ref{mainresult}, the solution corresponding to $q_1=q_2=r=1$ is
$$
u=-x_1-2x_2
$$
and the running cost is 
$$
L(x_1,x_2,u)=(x_1+x_2)^2+x_2^2+u^2
$$
The closed loop system is critically damped and has a double pole at $-1$. 
The value function is 
$$
V=\left(x_1+x_2\right)^2+x_2^2
$$
which can be rewritten as $V=x^TPx$ where
$$
P=\left[
\begin{array}{cc}
1 & 1\\
1 & 2
\end{array}
\right]
$$
Note that 
$$
\dot V=-\left(x_1+x_2\right)^2-x_2^2-\left(x_1+2x_2\right)^2 < 0,~\forall(x_1,x_2)\neq (0,0)
$$ 
Therefore, the value function is a global Lyapunov function.
\end{example}
\vspace{10pt}

\begin{example}
For system (\ref{nlsystemeq4}) with $a\neq 0,~d=ca^{-1}$, $q_1=q_2c^2a^{-2}$ and irrespectively of $f(x_2),~q_1,~q_2$ one has $k_1k_2^{-1}=-d,~k=0$ and the solution is a linear controller 
$$
u=-k_2(x_2-ca^{-1}x_1)
$$
The running cost and the value function are respectively
$$
L=rk_2^2(x_2-ca^{-1}x_1)^2+ru^2
$$ 
and 
$$
V=rb^{-1}k_2(x_2-ca^{-1}x_1)^2
$$
Note that in this case the two differential equations in (\ref{nlsystemeq4}) can be combined and the dynamics become $\dot z = bu$ where  $z=x_2-ca^{-1}x_1$. The controller is $u=-k_2z$, which makes the trajectories of $z$ converge to the origin exponentially.
In fact, this all makes sense because according to Theorem \ref{mainresult}, the trajectories are guaranteed to converge to the minimizers of $L$ given by the points in the set $\{(x_1,x_2):~x_2=ca^{-1}x_1\}$, for which the value function is zero. However, note that in this case the value function is not a Lyapunov function because $k=0$ and there is no guarantee that the trajectories converge to the origin. It is however a Lyapunov function for the dynamics of $z$. If $c=0$, which implies $d=0$, then $q_1=k_1=0$ and the trajectories will converge to the set of points $\{(x_1,x_2):~x_2=0\}$. But for $x_2=0$ we have $\dot x_1 = ax_1$ and $x_1$ therefore converges to zero if and only if $a<0$.  
\end{example}

\begin{example}
For system (\ref{nlsystemeq4}) consider $f(x_2)=x_2^3,~a=c=d=0$ and $b=1$. 
According to Theorem \ref{mainresult}, the optimal controller corresponding to $q_1=q_2=r=1$ is 
$$
u=-x_1-x_2-x_2^3
$$
the running cost is 
$$
L(x_1,x_2,u)=(x_1+x_2)^2+x_2^6+u^2
$$
and the value function is 
$$
V=\left(x_1+x_2\right)^2+0.5x_2^4
$$ 
Note that 
$$
\dot V=-\left(x_1+x_2\right)^2-x_2^6-\left(x_1+x_2+x_2^3\right)^2 < 0,~\forall(x_1,x_2)\neq (0,0)
$$
Therefore, the value function is a global Lyapunov function.
\end{example}
\vspace{10pt}

\begin{figure}[t] 
\centerline{ \epsfxsize=3.2in \epsffile{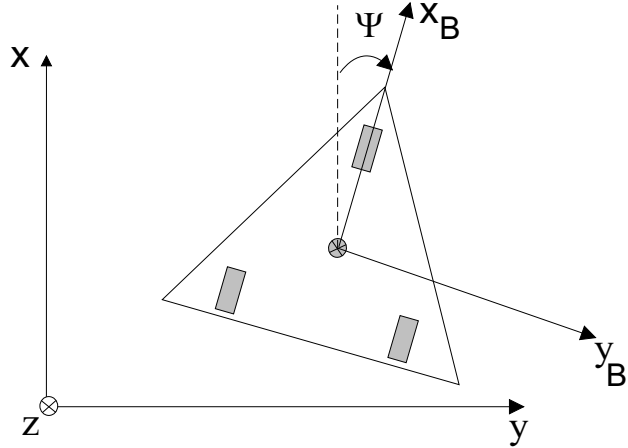} }
\caption{Path Following of Unicycle}
\label{pathfollowing}
\end{figure}

\begin{figure}[t] 
\centerline{ \epsfxsize=3.2in \epsffile{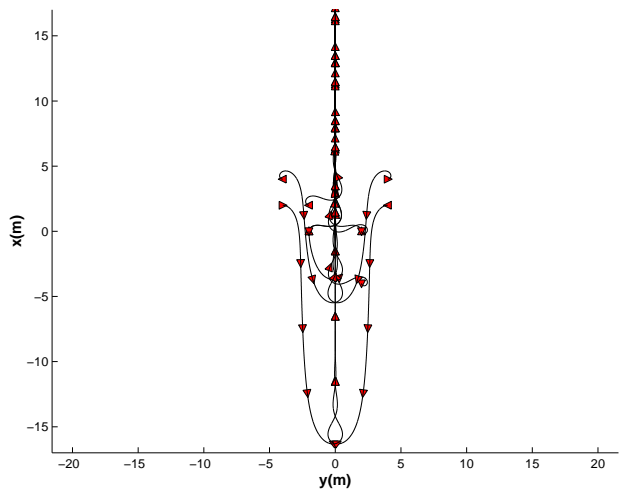} }
\caption{Unicycle Trajectories}
\label{unicycle_trajectories}
\end{figure}

\begin{example}\label{cart}
For system (\ref{nlsystemeq4}) let $f(x_2)=\sin(x_2),~a=c=d=0$ and $b=1$. 
This system is the kinematics model on the $x-y$ plane for path following of the line $y=0$ at constant unitary velocity by a unicycle. 
For this model, based on figure \ref{pathfollowing}, one has $x_1=y,~x_2=\psi$.
According to Theorem \ref{mainresult}, if $q_1=q_2=r=1$, the optimal controller is 
$$
u=-x_1-x_2-\sin(x_2)
$$ 
the running cost is 
$$
L(x_1,x_2,u)=(x_1+x_2)^2+\sin^2(x_2)+u^2
$$
and the value function is 
$$
V=\left(x_1+x_2\right)^2+2-2\cos(x_2)
$$ 
The derivative of the value function is 
$$
\begin{array}{rcl}
\dot V&=&-\left(x_1+x_2\right)^2-\sin^2(x_2)\\
&&-\left[x_1+x_2+\sin(x_2)\right]^2\le 0
\end{array}
$$
Therefore, the value function is a local Lyapunov function, which proves local stability in the sense of Lyapunov. However, the Lyapunov function is not radially unbounded (it is zero for $x_1=-x_2=2n\pi$ for $n$ integer) and asymptotic stability to the origin cannot be proved. In fact, by LaSalle's Invariance Principle \cite{khalil}, the trajectories are only guaranteed to converge to the largest invariant set contained in $\{(x_1,x_2):~\dot V(x_1,x_2)=0\}$, which is the set $\{(x_1,x_2):~x_1=-x_2,~x_2=n\pi\}$ where $n$ is an integer.
Notice that this is also the set of minimizers of $L$, which is in accordance with Theorem \ref{mainresult}.
Furthermore, invoking the result of Theorem \ref{mainresult}, one cannot guarantee convergence to the origin because $L$ is not convex in this case.
Figure \ref{unicycle_trajectories} shows several trajectories of the unicycle for different initial conditions. Convergence to the desired path is clearly seen for the initial conditions shown in the figure.
\end{example}

\section{Optimal Control Problem Definition and Solution: Third Order Systems}\label{threestates}
\noindent The results of the previous section are now extended to a class of third order systems for which $a=c=0$.
Consider the following optimal control problem
\beq
\begin{array}{rclr}
V(x_0)& = & \inf \int_0^{\infty} \left\{q_1x_1^2+q_2x_2^2+q_3x_3^2+Q(x)+ru^2\right\}dt&\\
& s.t.& \dot x_1(t) = f(x_2)&\label{nlsystemeq43d}\\
			&& \dot x_2(t) = df(x_2)+g(x_3)&\\
			&& \dot x_3(t) = bu&\\
			&& x(0)=x_0,~u\in{\mathcal U}&
\end{array}
\eeq
where it is assumed that $~q_1\ge 0,~q_2\ge0,~q_3>0,~r>0,~b\neq 0,~x(t)=[x_1(t)~x_2(t)~x_3(t)]^T\in\RE^3,~d\in\RE$.
The set $\mathcal U$ represents the allowable inputs, which are considered to be Lebesgue integrable functions. 
The functions $f,g$ are not identically zero and are assumed to be continuous with $f(0)=g(0)=0$.
The function $g(x_3)$ is assumed to have a locally positive definite anti-derivative $G(x_3)$ such that $G'(x_3)=g(x_3)$.

As before, we start by presenting necessary conditions that the value function $V$ must verify for a solution of the form (\ref{control3d}) to exist.
\newline
\begin{lemma}\label{necessary3d}
Assume that a control solution of the form 
\begin{equation}\label{control3d}
u(x)=-k_1x_1-k_2x_2-k_3x_3-k_4f(x_2)-k_5g(x_3)
\end{equation}
exists for problem (\ref{nlsystemeq43d}) and that a class ${\mathcal C}^1$ function $V$
exists that verifies the corresponding HJB equation
\beq\label{hjbvector3d}
\inf_u H\left(x_1,x_2,x_3,u,V_{x_1},V_{x_2},V_{x_3}\right)=0
\eeq
where
\beq\label{Hamiltonian3d}
\begin{array}{rcl}
H&=&q_1x_1^2+q_2x_2^2+q_3x_3^2+Q(x)+ru^2+V_{x_1}f(x_2)\\
&&+V_{x_2}g_3(x_2,x_3)+V_{x_3}bu
\end{array}
\eeq
with
$$
g_3=df(x_2)+g(x_3)
$$
and with boundary condition $V(0)=0$. 
Then $V$ must be of the form
\beq\label{valuefunction3d}
V(x)=2b^{-1}r\left(k_1x_1x_3+k_2x_2x_3+k_4f(x_2)x_3+k_3\frac{x_3^2}{2}+k_5G(x_3)\right)+h(x_1,x_2)
\eeq
where $h$ and $G$ are functions of class ${\mathcal C}^1$ with
\begin{eqnarray}
g(x_3) &=& G'(x_3)\label{antiderivative3d}\\
h(0,0) &=& -2b^{-1}rk_5G(0)\label{necessarycond1b3d}
\end{eqnarray}
\end{lemma}
{\bf Proof.}
Consider the HJB equation (\ref{hjbvector3d}) associated with (\ref{nlsystemeq43d}).
The necessary condition on $u$ to be a minimizer is
\beq\label{dvdx33d}
V_{x_3}=-2b^{-1}ru(x)
\eeq
and therefore
\beq\label{valuetwostates3d}
V(x)=-2b^{-1}r\int u(x)dx_3 + h(x_1,x_2)
\eeq
where $h$ is an arbitrary integration function of $x_1$ and $x_2$.
Searching for a solution of the form (\ref{control3d}), expression (\ref{valuetwostates3d}) becomes (\ref{valuefunction3d}).
From the boundary condition $V(0)=0$ one obtains the constraint (\ref{necessarycond1b3d}).
This finishes the proof. \done

\begin{theorem}\label{sufficient3d}
Let $Q(x)$ be chosen as
\beq\label{Q3d}
Q=2rk_1k_2x_1x_2+2rk_1k_3x_1x_3+2rk_2k_3x_2x_3+r\left(k_4^2-2dk_4k_5\right)f^2+rk_5^2g^2-2rb^{-1}k_4f'x_3\left(df+g\right)
\eeq
Then there exist gains $k_1,~k_2,~k_3,~k_4,~k_5$ verifying
\beq\label{k1constraints3d}
k_1=\sqrt{\frac{q_1}{r}}
\eeq
\beq\label{k2constraints3d}
k_2=\sqrt{\frac{q_2}{r}}
\eeq
\beq\label{k3constraints3d}
k_3=\sqrt{\frac{q_3}{r}}
\eeq
\beq\label{k4constraints3d}
k_4=b^{-1}k_3^{-1}\left(k_1+dk_2\right)
\eeq
\beq\label{k5constraints3d}
k_5=b^{-1}k_3^{-1}k_2
\eeq
such that the control input (\ref{control3d}) is a solution of the HJB equation (\ref{hjbvector3d}) associated with (\ref{nlsystemeq43d}) with value function
\begin{equation}\label{valuefunction3dfinal}
\begin{array}{rcl}
V(x)&=&rb^{-1}\left(\frac{k_1}{\sqrt{k_3}}x_1+\frac{k_2}{\sqrt{k_3}}x_2+\sqrt{k_3}x_3\right)^2\\
&&+2b^{-1}r\left[bk_4k_5\left(F(x_2)-F(0)\right)+k_4x_3f(x_2)+k_5\left(G(x_3)-G(0)\right)\right]
\end{array}
\end{equation}
where $u$ is given by (\ref{control3d}).
The function $V$ is also a local Lyapunov function for the system provided it is positive definite in a neighborhood of the origin and
\begin{equation}\label{Lpsd}
r\left(k_1x_1+k_2x_2+k_3x_3\right)^2+rb^{-2}k_3^{-2}\left(k_1^2-d^2k_2^2\right)f^2+rk_5^2g^2-2rb^{-1}k_4f'x_3\left(df+g\right)\ge 0
\end{equation}
If $V$ is globally positive definite and radially unbounded then it is a global Lyapunov function.
Finally, the trajectories converge to one of the minimizers of $L$. If $L$ is convex, then the trajectories will converge to the origin for all initial conditions.
\end{theorem}
{\bf Proof}.
Differentiating (\ref{valuefunction3d}) with respect to $x_1$ yields
\begin{equation}\label{dvdx13d}
V_{x_1}=2rb^{-1}k_1x_3+h_{x_1},
\end{equation}
and with respect to $x_2$ yields
\begin{equation}\label{dvdx23d}
V_{x_2}=2rb^{-1}\left(k_2+k_4f'\right)x_3+h_{x_2},
\end{equation}
where $h_{x_1}$ is the derivative of $h$ with respect to $x_1$ and $h_{x_2}$ is the derivative of $h$ with respect to $x_2$.
Replacing (\ref{control3d}), (\ref{dvdx13d}), (\ref{dvdx23d}), and (\ref{dvdx33d}) in (\ref{hjbvector3d}) yields after rearranging
\begin{equation}\label{hjbvector23d}
\begin{array}{rcl}
0&=&\left(q_1-rk_1^2\right)x_1^2+\left(q_2-rk_2^2\right)x_2^2+\left(q_3-rk_3^2\right)x_3^2+Q\\
&&-2rk_1k_2x_1x_2-2rk_1k_3x_1x_3-2rk_2k_3x_2x_3-rk_4^2f^2-rk_5^2g^2\\
&&-2r\left(k_4f+k_5g\right)\left(k_1x_1+k_2x_2+k_3x_3\right)-2rk_4k_5fg\\
&&+2rb^{-1}x_3\left[k_4f'\left(df+g\right)+f\left(k_1+dk_2\right)+gk_2\right]+\left(h_{x_1}+dh_{x_2}\right)f+h_{x_2}g
\end{array}
\end{equation}
Using (\ref{Q3d})--(\ref{k3constraints3d}), (\ref{hjbvector23d}) transforms to
\begin{equation}\label{hjbvector23d2}
\begin{array}{rcl}
0&=&-2r\left(k_4f+k_5g\right)\left(k_1x_1+k_2x_2+k_3x_3\right)-2rk_4k_5fg-2rdk_4k_5f^2\\
&&+2rb^{-1}x_3\left[f\left(k_1+dk_2\right)+gk_2\right]+\left(h_{x_1}+dh_{x_2}\right)f+h_{x_2}g
\end{array}
\end{equation}
Making
\begin{equation}\label{dhdx23d}
h_{x_2}=2rk_5\left(k_1x_1+k_2x_2\right)+2rk_4k_5f
\end{equation}
yields by integration
\begin{equation}\label{h3d}
h=2rk_5k_1x_1x_2+rk_5k_2x_2^2+2rk_4k_5F(x_2)+w(x_1)
\end{equation}
where $w$ is an arbitrary integration function of $x_1$.
Taking the derivative of (\ref{h3d}) with respect to $x_1$ yields
\begin{equation}\label{dhdx13d}
h_{x_1}=2rk_5k_1x_2+w'(x_1)
\end{equation}
Replacing (\ref{dhdx23d}) and (\ref{dhdx13d}) into (\ref{hjbvector23d2}), making
\begin{equation}\label{w3d}
w'(x_1)=2rk_1\left(k_4-dk_5\right)x_1
\end{equation}
and using (\ref{k4constraints3d})--(\ref{k5constraints3d}) yields the identity $0=0$ which proves that the HJB equation is satisfied.
This is a sufficient condition for the control input (\ref{control3d}) to be a solution that minimizes the cost of problem (\ref{nlsystemeq43d}) because the second derivative of the Hamiltonian (\ref{Hamiltonian3d}) with respect to $u$ is equal to $2r>0$.
Using (\ref{k1constraints3d})--(\ref{k5constraints3d}) the running cost is given by
\begin{equation}\label{running3d}
L=r\left(k_1x_1+k_2x_2+k_3x_3\right)^2+rb^{-2}k_3^{-2}\left(k_1^2-d^2k_2^2\right)f^2+rk_5^2g^2-2rb^{-1}k_4f'x_3\left(df+g\right)+ru^2\nonumber
\end{equation}
and it is non-negative with a minimun at $x_1=x_2=x_3=u=0$ under the assumption (\ref{Lpsd}).
Integrating (\ref{w3d}), using (\ref{h3d}), (\ref{k4constraints3d})--(\ref{k5constraints3d}) and the boundary condition $V(0)=0$, (\ref{valuefunction3d}) yields the value function (\ref{valuefunction3dfinal}).
Since $\dot V=-L\le 0$, the function $V$ is also a local Lyapunov function for the system if it is positive definite in a neighbourhood of the origin.
The rest of the proof follows the same argument as in the proof of Theorem \ref{case1}.\done
\begin{remark}
It is interesting to note the similarity in the strucure of (\ref{finalvalue}) and (\ref{valuefunction3dfinal}) for the case $c=0$. It is also worth to mention that for $d=0$ the results of Theorem \ref{sufficient3d} agree with the ones obtained in \cite{gholitabarrodrigues10}.
\end{remark}

%\section{Examples}\label{examples}
\begin{example}\label{wmr}
Consider now the third order integrator extension of example \ref{cart}.
The dynamics are  given by (\ref{nlsystemeq43d}) with $f(x_2)=\sin (x_2),~g(x_3)=x_3$, $b=1,~d=0$.
If $q_1 = q_3 = r = 1$ and $q_2 = 4$, then according to Theorem \ref{sufficient3d} the optimal controller is
\begin{equation}\label{eqC34}
u = -x_1 -2 x_2 -3 x_3 - sin(x_2)
\end{equation}
the running cost is
\begin{equation}\label{eqC35}
L = (x_1 + 2 x_2 + x_3)^2 + (4-2cos(x_2) )x_3^2+ \sin^{2}(x_2) + u^2
\end{equation}
and the value function is
\begin{equation}\label{eqC36}
V(x)=(x_1 + 2 x_2 + x_3)^2 + 2x_3^2 + 2x_3 \sin(x_2)- 4\cos(x_2)+4
\end{equation}
Computing the Hessian of $V$ and approximating the $\sin$ and $\cos$ functions by their first order Taylor series around zero one finds that $V$ is guaranteed to be positive definite in the set $\{(x_1,x_2,x_3)\in \mathbb{R}^3:~|x_2|<\epsilon,~|x_3|<1.5\epsilon^{-1}\}$ for small $\epsilon$.
Note that if one plots the functions that are the principal minors of the Hessian, one can actually find that $\epsilon$ can be as big as $\pi/10$ and the values of $x_3$ can still be obtained from the approximation above giving an accurate estimation of the region where $V$ is positive definite.
Moreover, the derivative of the value function is
\begin{equation}\label{eqC37}
\begin{split}
\dot V&=-(4-2cos(x_2))x_3^2 -(x_1 + 2 x_2 + x_3)^2 -  \sin^{2}(x_2) \\ 
     & - (x_1 +2 x_2 +3 x_3 + sin(x_2))^2 
\end{split}
\end{equation}
and is negative definite for $x_{2}\in (-\pi,\pi)$. Therefore, the value function is a local Lyapunov function in the largest invariant set contained in $\bigl \{ \thinspace (x_1,x_2,x_3) \in \mathbb{R}^3 \thinspace  | \thinspace |x_2| < \pi \bigr \}\cap\bigl \{ \thinspace (x_1,x_2,x_3) \in \mathbb{R}^3 | \thinspace V>0\}$ where $>0$ stands for positive definite. Note that, as in the previous example, one cannot guarantee convergence to the origin from any initial condition because $L$ is not convex.
\end{example}

\section{Conclusions}
This paper presented an inverse optimality method to solve a class of nonlinear optimal control problems. The method is inverse optimal because the running cost that renders the control input optimal is also explicitly determined.
The resulting running cost was shown to be a sensible non-negative cost with a minimum at the origin.

There are two main advantages of this method.
First, the analytical solution for the control input is obtained directly without needing to assume or compute a coordinate transformation, value function or Lyapunov function.
The value function and a Lyapunov function can however be computed after the control input has been found.
Another advantage is that it is capable of solving many examples of interest, inlcuding the Van der Pol oscillator, mass-spring systems and vehicle path following.
The main drawback of the method is that it is restricted to a specific class of optimal control problems for which the dynamics are affine in the input and the cost is quadratic in the input.

Two interesting conclusions can be drawn from this work.
First, the value function contains terms that are the negative integral of the control input. Regarding the control input as a force and the value function as potential energy, this integration leads to the usual expression for conservative forces, which is physically interesting.
Second, this work emphasizes the importance of cross terms on the state to find a solution to some optimal control problems. This is not only true in the value function, where they are needed to make the input be a function of both state variables, but also in the cost. Furthermore, making the cost depend on the nonlinearity, potentially including nonquadratic terms on the state, seems to be an important feature of this method.
This is in contrast to the traditional quadratic costs that have been used in a great percentage of the available literature in optimal control. 

\section*{Acknowledgments}
The authors would like to acknowledge the Natural Sciences and Engineering Research Council of Canada (NSERC) for funding this research and would like to dedicate this paper to the memory of Dr. El-Kebir Boukas.

\end{document}